\documentclass[10pt]{article}

\usepackage{amsmath,makeidx,amssymb,amsthm,amscd,index,mathset,psfig,graphicx}
\usepackage{color}
\usepackage[active]{srcltx}

\newcommand{\D}{\protect\displaystyle}
\newcommand{\T}{\protect\textstyle}
\newcommand{\eps}{\varepsilon}
\newcommand{\lbd}{\lambda}
\newcommand{\ipl}{\langle} 
\newcommand{\ipr}{\rangle}

\newcommand{\Lb}{\mathcal L}

\newcommand{\half}{{\textstyle\frac{1}{2}}}

\newtheorem{theorem}{Theorem}[section]

\newtheorem{lemma}[theorem]{Lemma}

\newtheorem{corollary}[theorem]{Corollary}
\newtheorem{remark}[theorem]{Remark}

\begin{document}
\setcounter{footnote}{2}

\title{\bf Regularization by dynamic programming}
\author{Stefan Kindermann\footnote{Johann
  Radon Institute for Computational and Applied Mathematics,
  Austrian Academy of Sciences, Altenbergerstrasse 69, A-4040 Linz,
  Austria. \text{email:  stefan.kindermann@oeaw.ac.at}}
   \and  A. Leit\~ao\footnote{Department of Mathematics, Federal University of
 St. Catarina, P.O. Box 476,  \mbox{88.040-900~Florianopolis}, Brazil. 
  \text{email: aleitao@mtm.ufsc.br} } }

\date{} \maketitle

\begin{abstract}
We investigate continuous regularization methods for linear inverse problems
of static and dynamic type. These methods are based on dynamic programming
approaches for linear quadratic optimal control problems.
We prove regularization properties and also obtain rates of convergence
for our methods.
A numerical example concerning a dynamical electrical impedance tomography
(EIT) problem is used to illustrate the theoretical results.
\end{abstract}

%---------------------------------------------------------------------------
% Section 1
%---------------------------------------------------------------------------
\section{Introduction} \label{sec:introd}

We begin by introducing the notion of dynamic inverse problems. Roughly
speaking, these are inverse problems in which the measuring process
--performed to obtain the data-- is time dependent. As usual, the problem
data corresponds to indirect information about an unknown parameter, which
has to be reconstructed. The desired parameter is allowed to be itself
time dependent.

Let $X$, $Y$ be Hilbert spaces. We consider the inverse problem of
finding $u: [0,T] \to X$ of the system
\begin{equation} \label{eq:dyn-ip}
F(t) u(t) \ = \ y(t)\, ,\ t \in [0,T]\, ,
\end{equation}
where $y: [0,T] \to Y$ are the dynamic measured data and $F(t): X \to Y$
are linear ill-posed operators indexed by the parameter $t \in [0,T]$.
Notice that $t \in [0,T]$ corresponds to a (continuous) temporal index.
The linear operators $F(t)$ map the unknown parameter $u(t)$ to the
measurements $y(t)$ at the time point $t$ during the finite time interval
$[0,T]$. We shall refer to (\ref{eq:dyn-ip}) as {\em dynamic inverse
problem}.

Since the operators $F(t)$ are ill-posed, at each time point $t \in [0,T]$
the solution $u(t)$ does not depend on a stable way on the right hand side
$y(t)$. Therefore, regularization techniques have to be used in order to
obtain a stable solution $u(t)$. In this article we consider time dependent
regularization methods \cite{KL06a}, which take into account the fact that
the parameter $u(t)$ evolves continuously with the time.

If the measuring process is stationary and the parameter is not time
dependent, the dynamic inverse problem (\ref{eq:dyn-ip}) reduces to the
standard problem of finding a solution $u \in X$ of the equation
\begin{equation} \label{eq:stat-ip}
F u \ = \ y\, ,
\end{equation}
where $F:X \to Y$ is a linear ill-posed parameter to output operator and
$y \in Y$. In opposition to (\ref{eq:dyn-ip}) we shall refer to
(\ref{eq:stat-ip}) as {\em static inverse problem}.

The second main goal in this article is to investigate continuous
regularization methods \cite{Tau94} for the inverse problem (\ref{eq:stat-ip}).
The regularization methods proposed for problems (\ref{eq:dyn-ip}) and
(\ref{eq:stat-ip}) are related by the fact that both of them derive
from a solution technique for linear quadratic optimal control problems
\cite{BK65}, the so-called {\em dynamic programming} \cite{Bel53,Bel57,BDS62}.

\subsubsection*{Some relevant applications}

As a first example of dynamic inverse problem, we present the {\em dynamical
source identification problem:} Let $u(x,t)$ be a solution to
$$ \Delta_x u(x,t) = f(x,t) \quad \mbox{ in } \Omega, $$
where $f(x,t)$ represents an unknown source which moves around and
might change shape with time $t$. The inverse problem in this case
is to reconstruct $f$ from single or multiple measurements of
Dirichlet and Neumann data $(u(x,t)$, $\partial_n u(x,t))$, on the
boundary $\partial\Omega$ over time $t \in [0,T]$.
Such problems arise in the field of medical imaging, e.g. brain
source reconstruction \cite{ABF02} or electrocardiography \cite{Le03}.

Many other 'classical' inverse problems have corresponding dynamic
counterparts, e.g., the {\em dynamic electrical impedance tomography problem}
consists in reconstructing the time-dependent diffusion coefficient
(impedance) in the equation
\begin{equation} \label{eq:dimp}
  \nabla_x \cdot(\sigma(.,t) \nabla_x) u(.,t) \ = \ 0,
\end{equation}
from measurements of the time-dependent Dirichlet to Neumann map
$\Lambda_\sigma$ (see the review paper \cite{CIN99}). This problem
can model a moving object with different impedance inside a fluid
with uniform impedance, for instance the heart inside the body.
Notice that in this case we assume the time-scale of the movement
to be large compared to the speed of the electro-magnetic waves.
Hence, the quasi-static formulation \eqref{eq:dimp} is a valid
approximation for the physical phenomena.

Another application concerning dynamical identification problems
for the heat equation is considered in \cite{KS97,KS99}.
Other examples of dynamic inverse problems can be found in
\cite{SLDBF01, SLWV02, SVVSK01, VK89, WB95}. In particular,
for applications related to process tomography, see the conference
papers by M.H.Pham, Y.Hua, N.B.Gray; M.Rychagov, S.Tereshchenko;
I.G.Kazantsev, I.Lemahieu in \cite{Pe00}.

\subsubsection*{Literature overview and outline of the paper}

Continuous and discrete regularization methods for static inverse problems
have been quite well studied in the last two decades and one can find
relevant information, e.g., in \cite{EHN96,EKN89,ES00,HNS95,Mor93,Tau94}
and in the references therein.

What concerns dynamic inverse problems, regularization methods were
considered for the first time in \cite{SL02, SLWV02}. There, the authors
analyze discrete dynamic inverse problems and propose a procedure called
{\em spatio temporal regularizer} (STR), which is based on the
minimization of the functional
\begin{equation} \label{eq:louis}
\Phi(u) \ := \ \T\sum\limits_{k=0}^{N} \| F_k u_k - y_k \|_{L^2}^2 +
     \lambda^2 \T\sum\limits_{k=0}^{N} \| u_k \|_{L^2}^2 +
         \mu^2 \T\sum\limits_{k=0}^{N-1} \frac{\| u_{k+1} - u_k \|_{L^2}^2}
                                      {(t_{k+1} - t_k)^2} .
\end{equation}
Notice that the term with factor $\lambda^2$ corresponds to the classical
(spacial) Tikho\-nov-Philips regularization, while the term with factor
$\mu^2$ enforces the temporal smoothness of $u_k$.

A characteristic of this approach is the fact that the hole solution
vector $\{ u_k \}_{k=0}^N$ has to be computed at a time. Therefore,
the corresponding system of equations to evaluate $\{ u_k \}$ has
very large dimension. In the STR regularization, the associated system
matrix is decomposed and rewritten into a Sylvester matrix form. The
efficiency of this approach is based on fast solvers for the Sylvester
equation.

In \cite{KL06a} continuous and iterative regularization methods based
on dynamic programming techniques were proposed as an alternative for
obtaining stable solutions of (\ref{eq:dyn-ip}). In this article the
authors verify regularization properties of the proposed methods and
present numerical realizations for a dynamic electrical impedance
tomography (EIT) problem, similar to the one treated in \cite{SLWV02}.

A word about the coupling of inverse problems and dynamic programming
theory. So far this theory have been mostly applied to solve particular
inverse problems. In \cite{KS97} the inverse problem of identifying the
initial condition in a semilinear parabolic equation is considered.
In \cite{KS99} the same authors consider a problem of parameter
identification for systems with distributed parameters.
In \cite{KL06b}, the dynamic programming methods are used in order
to formulate an abstract functional analytical method to treat static
inverse problems (\ref{eq:stat-ip}).

This paper is outlined as follows:
In Section~\ref{sec:deriv} we derive the solution methods discussed in
this paper.
In Section~\ref{sec:regul} we analyze some regularization properties of
the proposed methods. 
In Section~\ref{sec:appl} we present numerical realizations of the discrete
regularization method as well as a discretization of the continuous
regularization method. For comparison purposes we consider a dynamic EIT
problem, similar to the one treated in \cite{SLWV02}.

%---------------------------------------------------------------------------
% Section 2
%---------------------------------------------------------------------------
\section{Derivation of the regularization methods} \label{sec:deriv}

%---------------------------------------------------------------------------
\subsection{Static inverse problems} \label{ssec:s-ip}

We start this subsection defining an optimal control problem related
with the linear inverse problem (\ref{eq:stat-ip}). Let
$u_0 \in X$ be any approximation for the minimum norm solution
$u^\dag \in X$ of (\ref{eq:stat-ip}). We aim to find a function
$u: [0,T] \to X$ such that, $u(0) = u_0$ and
\begin{equation} \label{eq:ip-asymp}
\| F u(T) - y \| \ \approx \ \| F u^\dag - y \| \, .
\end{equation}
In the control literature, the function $u$ is called {\em trajectory}
(or {\em state}) and its evolution is is described by a dynamical system.
For simplicity, we choose a linear evolution model, i.e.
$u'  =  A u(t) + B v(t)$, $t \ge 0$,
where $A, B : X \to X$ are linear operators and $v:
[0,T] \to X$ is the {\em control} of the system.
Since our main concern is to satisfy the property in (\ref{eq:ip-asymp}),
it is enough for our purpose to consider a simpler dynamic, which does
not depend on the state $u$, but only on the control $v$ (for a dynamic including
the state see \cite{KiNa06}). This justifies
the choice of the dynamic: $u' = v$, $t \ge 0$. In this case, the control
$v$ corresponds to a {\em velocity function}.

The next step is to choose the objective function for our control problem.
The following choice is related to the minimization of both the residual
norm and the velocity norm along the trajectories
$$ J(u,v) := \T\frac{1}{2} \D\int_0^T \|F u(t) - y\|^2 + \| v(t) \|^2\ dt . $$
Putting all together we obtain the following abstract optimal control
problem in Hilbert spaces:
\begin{equation} \label{eq:ccp-hs}
  \left\{ \begin{array}{l}
      {\rm Mimimize} \ J(u,v) = \frac{1}{2} \displaystyle
                     \int_0^T \| F u(t) - y \|^2 + \| v(t) \|^2 \ dt \\
      {\rm s.t.} \\
      u' = v, \ t \ge 0 ,\ \ u(0) = u_0\, ,
   \end{array} \right.
\end{equation}
where the (fixed but arbitrary) final time $T > 0$ will play the role
of the regularization parameter. The functions $u,v : [0,T] \to X$
correspond respectively to the trajectory an the control of the system,
and the pairs $(u,v)$ are called {\em processes}.

Next we define the residual function $\eps(t) := F u(t) - y$ associated
to a given trajectory $u$. Notice that this residual function evolves
according to the dynamic
$$ \eps' = F u'(t) = F v(t)\, ,\ t \ge 0\, . $$
With this notation, problem (\ref{eq:ccp-hs}) can be rewritten in the
following form
\begin{equation} \label{eq:ccpr-hs}
  \left\{ \begin{array}{l}
      {\rm Mimimize} \ J(\eps,v) = \frac{1}{2} \D
                     \int_0^T \| \eps(t) \|^2 + \| v(t) \|^2 \ dt \\
      {\rm s.t.} \\
      \eps' = F v, \ t \ge 0 ,\ \ \eps(0) = F u_0 - y\, .
   \end{array} \right.
\end{equation}

In \cite[Proposition~2.1]{KL06b} the equivalence between the solvability
of the optimal control problem (\ref{eq:ccp-hs}) and the auxiliary problem
(\ref{eq:ccpr-hs}) is established.
In the sequel, we derive the dynamic programming approach for the
optimal control problem in (\ref{eq:ccpr-hs}). We start by introducing
the first Hamilton function. This is the function $H: \mathbb R \times
X^3 \to \mathbb R$ given by
$$  H(t,\eps,\lbd,v) \ := \ \D\ipl \lbd , F v \ipr +
    \T\frac{1}{2}\D [ \ipl \eps , \eps \ipr + \ipl v , v \ipr ] \, . $$
Notice that the variable $\lbd$ plays the role of a Lagrange multiplier
in the above definition. According to the Pontryagin's maximum principle,
the Hamilton function furnishes a necessary condition of optimality for
problem (\ref{eq:ccpr-hs}). Furthermore, since this function (in this
particular case) is convex in the control variable, this optimality
condition also happens to be sufficient. Recalling the maximum principle,
along an optimal trajectory we must have
\begin{equation} \label{eq:opt-cond}
0 \ = \ \frac{\partial H}{\partial v}(t,\eps(t),\lbd(t),v(t))
      \ = \ F^* \lbd(t) + v(t) \, .
\end{equation}
This means that the optimal control $\bar v$ can be obtained directly
from the Lagrange multiplier $\lbd: [0,T] \to X$, by the formula
$$  \bar v(t) = -F^* \lbd(t)\, ,\ \forall t\, . $$
Therefore, the key task is actually the evaluation of the Lagrange
multiplier. This leads us to the Hamilton-Jacobi equation. Substituting
the above expression for $\bar v$ in (\ref{eq:opt-cond}), we can define
the second Hamilton function $\mathcal H: \mathbb R \times X^2
\to \mathbb R$
$$  \mathcal H(t,u,\lbd) \ := \ \min_{v \in X}
    \{ H(t,\eps,\lbd,v) \} \ = \ \T\frac{1}{2} \D\ipl \eps , \eps \ipr -
    \T\frac{1}{2} \D\ipl \lbd , F F^* \lbd \ipr \, . $$
Now, let $V: [0,T] \times X \to \mathbb R$ be the value function
for problem (\ref{eq:ccpr-hs}), i.e.
\begin{eqnarray}
  V(t,\xi) \!\!\! & := & \!\!\!
                  \min\Big\{ \T\frac{1}{2} \int_t^T \| \eps(s) \|^2 +
                  \|v(s)\|^2 \, ds \ \Big| \ (\eps,v) \
                  {\rm admissible\ process} \nonumber \\
           \!\!\! &    & \!\!\!
                  {\rm \ \ \ \ \ \ \ \ for\ problem\
                  (\ref{eq:ccpr-hs})\ with\ initial\ condition} \
                  \eps(t) = \xi \Big\} \, . \label{def:cvf}
\end{eqnarray}
The interest in the value function follows
from the fact that this function is related to the Lagrange multiplier
$\lbd$ by the formula: $\lbd(t) = \partial V / \partial\eps (t,\bar \eps)$,
where $\bar \eps$ is an optimal trajectory.

From the control theory we know that the value function is a solution
of the Hamilton-Jacobi equation
\begin{equation} \label{eq:hjb}
 \frac{\partial V}{\partial t}(t,\eps) +
  \mathcal H(t,\eps,\frac{\partial V}{\partial \eps}(t,\eps)) \ = \ 0\, .
\end{equation}
Now, making the ansatz: $V(t,\eps) = \frac{1}{2} \ipl \eps , Q(t) \eps \ipr$,
with $Q: [0,T] \to \mathbb R$, we are able to rewrite (\ref{eq:hjb}) in
the form
$$  \ipl \eps , Q'(t) \eps \ipr + \ipl \eps , \eps \ipr - 
    \ipl Q(t) \eps , F F^* Q(t) \eps \ipr \ = \ 0 \, . $$
Since this equation must hold for all $\eps \in X$, the function
$Q$ can be obtained by solving the Riccati equation
\begin{equation} \label{eq:riccati}
Q'(t) \ = \ -I + Q(t) F F^* Q(t) \, .
\end{equation}
Notice that the cost of all admissible processes for an initial
condition of the type $(T,\eps)$ is zero. Therefore we have to
consider the Riccati equation (\ref{eq:riccati}) with the final
condition
\begin{equation} \label{eq:riccati-ic}
Q(T) \ = \ 0\, .
\end{equation}

Once we have solved the initial value problem (\ref{eq:riccati}),
(\ref{eq:riccati-ic}), the Lagrange multiplier is given by $\lbd(t)
= Q(t) \bar\eps(t)$ and the optimal control is obtained by the
formula $\bar v(t) = - F^* Q(t) \bar\eps(t)$. Therefore, the optimal
trajectory of problem (\ref{eq:ccp-hs}) is defined via
\begin{equation} \label{eq:cont-reg}
 \bar u' = - F^* Q(t) [F \bar u(t) - y] \, ,\ \ \bar u(0) = u_0 \, .
\end{equation}

We use the optimal trajectory defined by the initial value problem 
(\ref{eq:cont-reg}) in order to define a family of reconstruction
operators $R_T: X \to X$, $T \in \mathbb R^+$,
\begin{equation} \label{eq:cont-Rt}
R_T(y) \ := \ \bar u(T) \ = \ u_0 -
              \int_0^T F^* Q(t) [F \bar u(t) - y] \ dt\, .
\end{equation}
We shall return to the operators $\{ R_T \}$ in Section~\ref{sec:regul}
and prove that the family of operators defined in (\ref{eq:cont-Rt})
is a regularization method for (\ref{eq:stat-ip}) (see, e.g., \cite{EHN96}).

%---------------------------------------------------------------------------
\subsection{Dynamic inverse problems} \label{ssec:d-ip}

In the sequel we consider the dynamic inverse problem described in
(\ref{eq:dyn-ip}). As in the previous subsection, we shall look for
a continuous regularization strategy.

We start by considering the constrained optimization problem
\begin{equation} \label{eq:ccp}
  \left\{ \begin{array}{l}
      {\rm Mimimize} \ J(u,v) := \half \int_0^T \big[ \,
                     \ipl F(t)u(t)-y(t), \ L(t) [F(t)u(t)-y(t)] \ipr \\[1ex]
      \hskip4.2cm      + \ \ipl v(t), M(t) v(t) \ipr \, \big] \ dt \\[1ex]
      {\rm s.t.} \ \ u' = A(t) u + B(t) v(t), \ t \in [0,T] ,\ u(0) = u_0\, ,
   \end{array} \right.
\end{equation}
where $F(t)$, $u(t)$ and $y(t)$ are defined as in (\ref{eq:dyn-ip}),
$v(t) \in X$, $t \in [0,T]$, $L(t): Y \to Y$, $M(t): X \to X$,
$A(t) \equiv I: X \to X$, $B(t) \equiv 0$ and $u_0 \in X$.

Following the footsteps of the previous subsection, we define the first
Hamilton function $H: [0,T] \times X^3 \to \mathbb R$ by
$$  H(t,u,\lbd,v) \ := \ \ipl \lbd , v \ipr + \half
    [ \ipl F(t) u - y(t), L(t) (F(t) u - y(t)) \ipr + \ipl v, M(t)v \ipr ]
    \, . $$
Thus, it follows from the maximum principle: $0 = \partial H / \partial v
(t,u(t),\lbd(t),v(t)) = \lbd(t) + M(t) v(t)$, and we obtain a relation
between the optimal control and the Lagrange parameter, namely:
$\bar v(t) = - M^{-1}(t) \lbd(t)$.

As before, we define the second Hamilton function $\mathcal H: \mathbb R
\times X^2 \to \mathbb R$
$$  \mathcal H(t,u,\lbd) \, := \,
%   \min_{v \in X} \{ H(t,u,\lbd,v) \} \, = \,
    \half \ipl F(t) u - y(t) , L(t) (F(t) u - y(t) \ipr -
    \half \ipl \lbd , M(t)^{-1} \lbd \ipr \, . $$
Since $\lbd(t) = \partial V / \partial u(t,u)$, where $V: [0,T] \times X
\to \mathbb R$ is the value function of problem (\ref{eq:ccp}), it is
enough to obtain $V$. This is done by solving the Hamilton--Jacobi (HJ)
equation (see \eqref{eq:hjb})
\begin{eqnarray*}
%0 & = & V_t(t,u) + \mathcal H(t, u, V_u(t,u)) \\
 0 & = & V_t + \half \ipl F(t)u - y(t) , L(t) (F(t)u - y(t)) \ipr
             - \half \ipl V_u , M(t)^{-1} V_u \ipr \, .
\end{eqnarray*}
We make the ansatz $V(t,u) =
\half \ipl u , Q(t) u \ipr + \ipl b(t), u \ipr + g(t)$,
with $Q: [0,T] \to \mathbb R$, $b: [0,T] \to X$ and $g: \mathbb R \to
\mathbb R$. Then, we are able to rewrite the HJ equation above in the
form of a polynomial equation in $u$. Moreover, the quadratic, the
linear and the constant terms of this polynomial equation must all
vanish. Thus we obtain
\begin{equation} \label{eq:evol}
Q' = Q^* M(t)^{-1} Q - F^*(t) L(t) F(t) , \ \ 
b' = Q(t)^* M(t)^{-1} b + F^*(t) L(t) y(t) .
\end{equation}
%
% $$   g' = \half \ipl b(t), M(t)^{-1} b(t) \ipr
%         + \ipl y(t) , y(t) \ipr \, . $$
%
The final conditions $Q(T) = 0$, $b(T) = 0$ are derived just like in the
previous subsection.%
\footnote{Since function $g$ is not needed for the computation of
the optimal trajectory, we omit the expression of the corresponding dynamic.}
Once the above system is solved, the optimal control $\bar u$ is obtained
by solving
\begin{equation} \label{eq:u}
\bar u'(t) = -M^{-1}(t) V_u(t,u) = -M^{-1}(t) [ Q(t) \bar u(t) + b(t)]
\end{equation}
with initial condition $\bar u(0) = u_0$.

Following the ideas of the previous tutorial subsection, we shall
choose a family of operators $\{ M_\alpha, L_\alpha \}_{\alpha>0}$ and
use the corresponding optimal trajectories $\bar u_\alpha$ in order to
define a family of reconstruction operators $R_\alpha: L^2((0,T) ; Y)
\to H^1( (0,T) ; X)$,
$$ R_\alpha(y) := u_0 - \T\int_0^t
                  M_\alpha^{-1}(s) [ Q(s) \bar u(s) + b(s)] \, ds \, . $$

The regularization properties of the operators $\{ R_\alpha \}$ will be
analyzed in Section~\ref{sec:regul}.

%---------------------------------------------------------------------------
% Section 3
%---------------------------------------------------------------------------
\section{Analysis of regularization properties} \label{sec:regul}

%---------------------------------------------------------------------------
\subsection{Static inverse problems} \label{ssec:s-reg}

In this section we investigate the regularization properties
of the operator $R_T$ introduced in (\ref{eq:cont-Rt}).
Consider the Riccati equation (\ref{eq:riccati}) for the operator $Q$: 
We may express the operator $Q(t)$ via the spectral
family $\{F_\lambda\}$ of $FF^*$
(see e.g. \cite{EHN96}).
Hence, we make the ansatz
$$ Q(t) = \int q(t,\lambda) d F_\lambda \; . $$
Assuming that $q(t,\lambda)$ is $C^1$ we may find
from (\ref{eq:riccati}) together with the boundary condition at $t = T$
that
$$  \int \left( \T\frac{d}{d t} q(t,\lambda)  
    + 1 -  q(t,\lambda)^2 \lambda  \right) \ d F_\lambda  = 0, \ \ \ 
    q(T,\lambda) = 0 . $$
Hence, we obtain an ordinary differential equation  for $q$:
\begin{equation}\label{ode}
 \T\frac{d}{d t} q(t,\lambda) = -1 + \lambda q(t,\lambda)^2
\end{equation}
The solution to these equations is given by
\begin{equation}\label{qfunc}
 q(t,\lambda) = - \frac{1}{\sqrt{\lambda}} 
\tanh(\sqrt{\lambda}(t -T)) = 
 \frac{1}{\sqrt{\lambda}} 
\tanh(\sqrt{\lambda}(T-t)) .
\end{equation}
If $t < T$, then $Q(t)$ is nonsingular, since
$\lim_{x\to 0 } \frac{\tanh(x a)}{x}  = a$ and 
$\frac{\tanh( a x)}{x} $ is monotonically decreasing 
for $x >0$. Hence the spectrum of $Q(t)$ is
contained in the interval $[\frac{\tanh((T-t)\|F\|)}{\|F\|}, (T-t)]$.
Now consider the evolution equation (\ref{eq:cont-reg}):
The operator $Q(t)$ can be expressed as
$Q(t) = q(t,FF^*)$;
by usual spectral theoretic properties  (see, e.g., \cite{EHN96})
it holds that
$$ F^*q(t,F F^*)  =  q(t,F^* F)F^*. $$
Hence we obtain the problem
\begin{eqnarray}
 u'(t) &=&  -q(t,F^*F) \left( F^*F u(t) - F^*y \right)  \label{ip1}\\
 u(0) &=& u_0 \label{ip2}
\end{eqnarray}
%

% By linearity we may write  solution to the initial value problem
% (\ref{ip1},\ref{ip2}) 
% $u = u_h + u_i$,
% where $u_h$ solves (\ref{ip1}) with the initial condition
% $u_h(0) = 0$ and  
% $u_i$ solves
% \begin{eqnarray}\label{initial}
% u'(t) &=&  -q(t,F^*F) F^*F u(t)  \\
% u(0) & = & u_0 
% \end{eqnarray}

We may again use an ansatz via spectral calculus: if we set
$$ u(t) = \int g(t,\lambda) d E_\lambda F^* y$$
where $E_\lambda$ is the spectral family of $F^* F$,
we derive an ordinary differential equation for $g$.
%
%$$ \int_{\sigma} \frac{d}{dt} g(t,\lambda) 
%d E_\lambda F^*y = 
%\int_{\sigma} - q(t,\lambda) \lambda g(t,\lambda) 
%+ q(t,\lambda) d E_\lambda F^*y $$
%Hence
%$$ g'(t,\lambda) = - q(t,\lambda) \lambda g(t,\lambda) + q(t,\lambda)$$
%with $q$ as in (\ref{qfunc})
%
Similar as above, we can express the solution to (\ref{ip1},\ref{ip2}) 
in the form
%
%$g(0,\lambda) = 0$  is given by
%\begin{equation}
% g(t,\lambda) =  \frac{1 - \frac{\cosh(\sqrt{\lambda}(T-t))}{
%\cosh(\sqrt{\lambda} T )}}{\lambda}
%\end{equation}
%Similar we find for $u_i$,
%$$ u_i = \int_{\sigma} h(t,\lambda) d E_\lambda u_0$$
%with 
%$$ h(t,\lambda) = \frac{\cosh(\sqrt{\lambda}(T-t))}{\cosh(\sqrt{\lambda}T )}
%$$
%
\begin{equation}\label{tsol}
u(t) = 
\int \frac{1 - \frac{\cosh(\sqrt{\lambda}(T-t))}{
\cosh(\sqrt{\lambda} T)}}{\lambda} d E_\lambda F^*y 
+ \int
 \frac{\cosh(\sqrt{\lambda}(T-t))}{\cosh(\sqrt{\lambda} T )} d E_\lambda u_0.
\end{equation}
Setting $t=T$ we find an approximation of the solution 
\begin{equation}\label{solution}
u_T:=  u(T) = 
\int \frac{1 - \frac{1}{
\cosh(\sqrt{\lambda} T )}}{\lambda} d E_\lambda F^*y 
+ \int
 \frac{1}{\cosh(\sqrt{\lambda} T )} d E_\lambda u_0 .
\end{equation}
Note the similarity to Showalter`s methods \cite{EHN96},
where the term $\exp(\lambda T)$ instead of
$\cosh( \sqrt{\lambda} T )$ appears.

\begin{theorem} \label{th:convcont}
The operator $R_T$ in (\ref{eq:cont-Rt}) is a 
regularization operator with qualification $\mu_0 = \infty$:

If the data are exact, $y = F u^\dagger$ and 
$u^\dagger$ satisfies a source condition for some
$\mu > 0$ 
\begin{equation}\label{sourcecond} 
 \exists\ \omega \in X : \quad  u^\dagger  = (F^*F)^\mu \omega ,
\end{equation}
we have the estimate
$$ \|u_T - u^\dagger\| \leq C_\mu T^{-2 \mu} $$

If the data are contaminated with noise,
$ \| y -y_\delta\| \leq \delta $ and $y = F u^\dagger$ 
with $u^\dagger$ as in  (\ref{sourcecond}), then we have
$$ \|u_{T,\delta} - u^\dagger\| \leq
C_\mu T^{-2 \mu} + \delta T . $$
In particular, the a-priori parameter choice
$T \sim \delta^{\frac{-1}{2 \mu +1}} $ yields the optimal
order convergence rate
$$ \|u_{T,\delta} - u^\dagger\| \sim \delta^{\frac{2 }{2 \mu +1}} . $$
\end{theorem}

\noindent {\em Proof:} See \cite[Theorem~3.1]{KL06b}.
\bigskip

Comparing the dynamic programming approach with the Showalter method,
one observes that they are quite similar, with $T_{dyn}^2 \sim T_{Sw}$.
Hence, to obtain the same order of convergence we only need $\sqrt{T_{Sw}}$
of the time for the Showalter method.

%---------------------------------------------------------------------------
\subsection{Dynamic inverse problems} \label{ssec:d-reg}

Before we examine the regularization properties of the method derived
in Subsection~\ref{ssec:d-ip}, we state a result about existence and
uniqueness of the Riccati equations~\eqref{eq:evol}.

\begin{theorem}
If $F$, $L$, $M \in C([0,T],\Lb(X,Y))$, then the Riccati equation
\eqref{eq:evol} has a unique symmetric positive semidefinite solution
in $C^1([0,T],\Lb(X))$.
\end{theorem}

\noindent {\it Proof:} See \cite[Theorem~3.1]{KL06a}.

\begin{remark}
It is well known in control theory that the existence of a solution to
\eqref{eq:evol} can be constructed from the functional
\begin{multline} \label{eq:value}
V(t,\xi) := % \min_{{u(t) = \xi} \atop {u \in H^1([t,T],X)}}
            \min_{ \substack{ {u(t) = \xi} \\ {u \in H^1([t,T],X)} } }
            \half
            \T\int_t^T \ipl F(s) u(s) - y(s), L(s)[F(s) u(s) - y(s)] \ipr
            \\ + \ipl u'(s), M(s) u'(s) \ipr ds.
\end{multline}
This functional is quadratic in $u$ and, from the Tikhonov regularization
theory (see, e.g., \cite{EHN96}), it admits a unique solution $u$, and is
quadratic in $\xi$. Furthermore, the leading quadratic part $(\xi , Q(t)\xi)$
is a solution to the Riccati Equation.
\end{remark}

Next we consider regularization properties of the method derived in
Subsection~\ref{ssec:d-ip}. The following lemma shows that the solution
$u$ of \eqref{eq:u} satisfies the necessary optimality condition for the
functional
\begin{equation} \label{eq:cost-func}
 J(u) \ = \ \half
            \T\int_0^T \ipl F(s) u(s) - y(s), L(s)[F(s) u(s) - y(s)] \ipr
            \\ + \ipl u'(s), M(s) u'(s) \ipr ds
\end{equation}
(notice that this is the cost functional $J(u,v)$ in \eqref{eq:ccp} with
$v = u'$).

\begin{lemma}
Let $Q(t)$, $b(t)$, $u(t)$ be defined by \eqref{eq:evol}, \eqref{eq:u},
together with the boundary conditions $Q(T) = 0$, $b(T) = 0$ and
$u(0) = u_0$. Then, $u(t)$ solves
\begin{equation} \label{eq:diff}
F^*(t) L(t) F(t) u(t) - \left(M(t) u(t)'\right)' = F^*(t) L(t) y(t),
\end{equation}
together with the boundary conditions $u(0) = u_0$, $u'(T) = 0$.
\end{lemma}
\noindent {\it Proof:} See \cite[Lemma~3.3]{KL06a}.
\bigskip

Since the cost functional in \eqref{eq:cost-func} is quadratic, the
necessary first order conditions are also sufficient. Thus, the solution
$u(t)$ of \eqref{eq:diff} is actually a minimizer of this functional.
Including the boundary conditions we obtain the following corollary:

\begin{corollary} \label{cor:cont-min}
The solution $u(t)$ of \eqref{eq:diff} is a minimizer of the Tikhonov
functional in \eqref{eq:cost-func} over the linear manifold
$$ \mathcal H:= \{ u \in H^1([0,T],X) \ | \ u(0) = u_0 \} \, . $$
\end{corollary}

In particular, this means that the above procedure is a regularization
method for the inverse problem \eqref{eq:dyn-ip}.
Bellow we summarize a stability and convergence result. The proof uses
classical techniques from the analysis of Tikhonov type regularization
methods (cf. \cite{EHN96}, \cite{EKN89}) and thus is omitted.

\begin{theorem}
Let $M(t) \equiv \alpha I$, $\alpha > 0$, $L(t) > 0$, $t \in [0,T]$
and $J_\alpha$ be the corresponding Tikhonov functional given by
\eqref{eq:cost-func}. \\
\underline{Stability:} Let the data $y(t)$ be noise free and denote by
$u_\alpha(t)$ the minimizer of $J_\alpha$.
Then, for every sequence $\{\alpha_k\}_{k \in \mathbb{N}}$ converging to zero,
there exists a subsequence $\{\alpha_{k_j}\}_{j \in \mathbb{N}}$, such that
$\{ u_{\alpha_{k_j}} \}_{j \in \mathbb{N}}$ is strongly convergent. Moreover,
the limit is a minimal norm solution. \\
\underline{Convergence:} Let $\|y^\delta(t) - y(t)\| \leq \delta$.
If $\alpha = \alpha(\delta)$ satisfies
$$ \lim_{\delta \to 0} \alpha (\delta) = 0 \mbox{ \ \ and \ \ }
   \lim_{\delta \to 0} \delta^2 / \alpha (\delta) = 0 \, . $$
Then, for a sequence $\{\delta_k\}_{k\in\mathbb{N}}$ converging to zero, there
exists a sequence $\{\alpha_k := \alpha(\delta_k)\}_{k \in \mathbb{N}}$
such that $u_{\alpha_k}$ converges to a minimal norm solution.
\end{theorem}

For the sake of completeness we also include the corresponding algorithms for
the discretized dynamical case. Instead of having a continuous time
variable we assume that the data and the operator are given on 
discrete time steps $k=1,\ldots N$. Instead of \eqref{eq:cost-func} 
a functional is used where the integrals are replaced by sums, the 
derivatives by differences very similar to \eqref{eq:louis}. The operators
$F(t)$, $L(t)$ are replaced by sequences $F_k$, $L_k$, $k=1, \ldots N$, the
data $y(t)$ are now given as a sequence $y_k$ and instead of a time dependent solution
$u(t)$ we are looking for a set of solutions $u_k$. The dynamic programming
approach goes through in a similar manner and as a result we obtain an
iterative procedure instead of the evolution equations  \eqref{eq:evol},  \eqref{eq:u}.
The details can be found in \cite{KL06a,KL06b}, we just state the iterations. 
Set $Q_{N+1}:=0$, $b_{N+1}:= 0$ and fix an initial guess $u_0$. The corresponding
minimizer of the discrete version of \eqref{eq:cost-func} can be found by two
backwards iterations on $Q$ and $b$ and a forward iteration on $u$ (for simplicity
we put $M = I$):
% \[
% \begin{array}{lcrr}
% Q_{k-1} & = & (Q_k + \alpha^{-1} I)^{-1} Q_k + F_{k-1}^* L_{k-1} F_{k-1}
%               \qquad & k = N +1 \ldots 2 \\
% b_{k-1} & = &  (Q_k + \alpha^{-1} I)^{-1} b_k - F_{k-1}^* L_{k-1} y_{k-1}
%               \qquad & k = N+1 \ldots 2  \\
% u_{k} & = & (Q_{k} + \alpha I)^{-1} (\alpha u_{k-1} - b_{k})
%               \qquad & k = 1 \ldots N 
% \end{array}
% \]
\begin{align}
 Q_{k-1} & =  (Q_k + \alpha^{-1} I)^{-1} Q_k + F_{k-1}^* L_{k-1} F_{k-1}&
               \quad   & k  = N +1, \ldots, 2 \label{itone}\\
 b_{k-1} & =   (Q_k + \alpha^{-1} I)^{-1} b_k - F_{k-1}^* L_{k-1} y_{k-1} &
               \quad   & k = N+1, \ldots, 2  \label{ittwo} \\
 u_{k} & =  (Q_{k} + \alpha I)^{-1} (\alpha u_{k-1} - b_{k}) &
                  & k = 1, \ldots, N  \label{itthree}
\end{align}
It can be shown that the sequence $u_k$ defined in this way satisfies the
optimality conditions for \eqref{eq:cost-func} and  that
it is a minimizer similar as in Corollary~\ref{cor:cont-min}. By standard
Tikhonov theory this implies that this iterative procedure is a regularization.

%---------------------------------------------------------------------------
% Section 4
%---------------------------------------------------------------------------
\section{Application to a dynamical EIT problem} \label{sec:appl}

%---------------------------------------------------------------------------
\subsection{The model}
As a motivation for considering dynamical inverse problem we 
stated the dynamical impedance tomography problem, namely
to identify a time-dependent conductivity coefficient 
in \eqref{eq:dimp}, from the measurements of the
associated Dirichlet-to-Neumann (DN) map $\Lambda_\sigma$.
To be more concrete we assume that for a fixed time $t$
$u(.,t)$ is a solution to \eqref{eq:dimp}  on a
fixed domain $\Omega$, with Dirichlet data $f \in H^{\frac{1}{2}}
(\partial\Omega)$.
For any $f$ we can measure the associated Neumann data
$g = \frac{\partial}{\partial n} u(.,t) \in H^{-\frac{1}{2}}
(\partial\Omega)$. The knowledge of all pairs of Cauchy data $(f,g)$ is
equivalent to knowing the DN map $\Lambda_{\sigma(.,t)}: f \to g$.
Note that the equation  \eqref{eq:dimp} does not involve derivatives of $t$
and the time-dependence of $u$ and hence $\Lambda_{\sigma(.,t)}$ is only
introduced by the time-dependence of the coefficient $\sigma(.,t)$.

The inverse problem associated to dynamical EIT is to identify the parameter
$\sigma(.,t)$ on $\Omega \times [0,T]$ from the time-dependent DN
map $\Lambda_{\sigma(.,t)}$. Hence, in our notation to parameter-to-data map
$F(\sigma)$ is $\sigma \to \Lambda_{\sigma(.,t)}$. However, this map is
nonlinear and does not fit into the framework of our work, which right now
only deals with linear operators. It is therefore necessary to linearize the
problem. We assume that the conductivity coefficient is a small
(time-dependent) perturbation of a constant background conductivity:
$\sigma(x,t) = 1+ \gamma(x,t)$. In this case it makes sense to subtract the
constant-conductivity  operator $\Lambda_1$ from the data and linearize the
parameter-to-data map:
\[ \gamma(x,t) \to \Lambda_{1+\gamma(.,t)} - \Lambda_1 \sim F'(1) \gamma, \]
Here $F'(1)$ denotes the Fr\'echet-derivative  of the nonlinear
parameter-to-data map at conductivity $1$. Now the forward operator $F(t)$ in
\eqref{eq:dyn-ip} can be identified with $F'(1)$ and the problem fits into the
framework of the linear dynamical inverse problems. Note that in this case the
forward operator does not depend on time, but the data and the solution do, so
that we have a problem of the form $F u(t) = y(t)$, which of course can be
handled by the dynamic programming approach. Since we use simulated data
$y(t)$ we can either consider linearized data perturbed with random noise
$y(t) = F'(1) \gamma(t) + \mbox{noise}$, or we can as well take the nonlinear
data $y(t) = \Lambda_{1+\gamma(.,t)} - \Lambda_1$ and treat the linearization
error as a data error.

The linearized operator DN operator $F'(1) \gamma$ has the following form:
It maps the Dirichlet values $f$ to the Neumann values $g= \frac{\partial}
{\partial n} w |_{\partial \Omega}$, with $w$ the solution of of the
linearized problem of the form
\begin{eqnarray}
\Delta   w & = & - \nabla\cdot(\gamma(.,t) \nabla  u_0)
                   \quad w|_{\partial \Omega} = 0 \label{eq:dimp_lin} \\
\Delta u_0 & = & 0 \quad u_0|_{\partial \Omega} = f. \nonumber
\end{eqnarray}

\begin{remark}
Let us note that the dynamical inverse problems approach can also be used as
a {\em dimension reduction}, in the sense that we can solve a
three-dimensional problem by considering it as a two-dimensional one with
a parameter dependent operator and solution. This parameter represents the
third dimension.
If the dependence of the forward operator on the third dimension is low we
can interpret the third dimension as a time-variable and use the framework
in this paper to solve a three-dimensional problem by using only the
corresponding operators for the planar case, which obviously is much
simpler. Such an approach works, if the the three-dimensional problem can
be approximated by two-dimensional 'slices', as it is usually done in
computerized tomography or electrical impedance tomography.

But contrary to the standard approach where each 2D-slice is treated
separately, our approach allows a coupling of the solutions in each slice, to
build a continuous 3D solution. Moreover, we do not have the restriction that
each two-dimensional problem is of same type, since we allow our operator $F$
to depend on time (or the third space dimension). This can happen if, for
instance, the geometry of the problem changes with the third dimension.
All in all we think that the dynamical programming framework can improve
a 2D-slicing approach.
\end{remark}

Let us consider the computational aspect of the problem.
For the numerical setting it is more convenient to work with the
Neumann-to-Dirichlet (ND) operator instead of the DN, since it is a smoothing
operator, whereas the latter is not.

For a discretization we use piecewise linear finite element $\phi_i$ for
the solutions of the differential equation \eqref{eq:dimp} or its linearized
version \eqref{eq:dimp_lin}. For $\gamma$ we use piecewise constant
elements. The discretization of the boundary function is simply obtained by
using the boundary trace of the finite elements $\tilde{\phi_i} =
\phi_i|_{\partial \Omega}$.

It is well know that with finite elements the Neumann problem has a discrete
form
$$ \left( \begin{array}{cc} A_{11} & A_{12} \\
 A_{21} & A_{22} \end{array} \right) \left(\begin{array}{c}
  u_i \\ u_b \end{array} \right)  =
  \left(\begin{array}{c}
  0 \\ M g \end{array}\right), $$
where the matrices $A_{11},A_{12},A_{21},A_{22}$ are sub-matrices of the
stiffness matrix $A_{i,j} = \int_\Omega \gamma \nabla \phi_i \nabla \phi_j$
with  respect of a splitting of the indices into the interior and boundary
components. The matrix $M$ is coming from the contribution of the
Neumann-data $ g = \sum g_i \tilde{\phi}_i $ in the discretized  equations and
is defined as 
\begin{equation} \label{bound}
 M_{i,l} = \T\int_{\partial \Omega} \tilde{\phi_l} \tilde{\phi}_i d\sigma.  
\end{equation}
It is easy to show that the discretized ND map  $G$ is given 
the Schur-Complement of the stiffness matrix with respect to the interior
components:
\begin{equation} \label{Gdef}
 G := \left( A_{22} - A_{21} A_{11}^{-1} A_{12} \right)^{-1} M. 
\end{equation}
Since $G$ is in the data space we need a Hilbert space to measure 
the error. For this task the Hilbert-Schmidt inner product can be used as
in \cite{KL06a}. If $G_1,G_2$ are discretized ND maps we 
use as inner product \[ (G_1,G_2) = \mbox{trace}(G1\,.G2). \]

\subsection{Numerical results}

In this section we present the numerical results for the dynamical impedance
tomography.
We consider Equation~\eqref{eq:dimp} on the unit ball with a time-dependent
conductivity parameter $\sigma(x,t) = 1 + \gamma(x,t)$, where $\gamma$
represents an inhomogeneity. For the first example we used noise-free data
for the linearized problem, i.e. $y(.,t) = F'(1) \gamma(.,t)$.
 
\begin{figure}[ht]
\begin{center}
  \includegraphics[width=\textwidth]{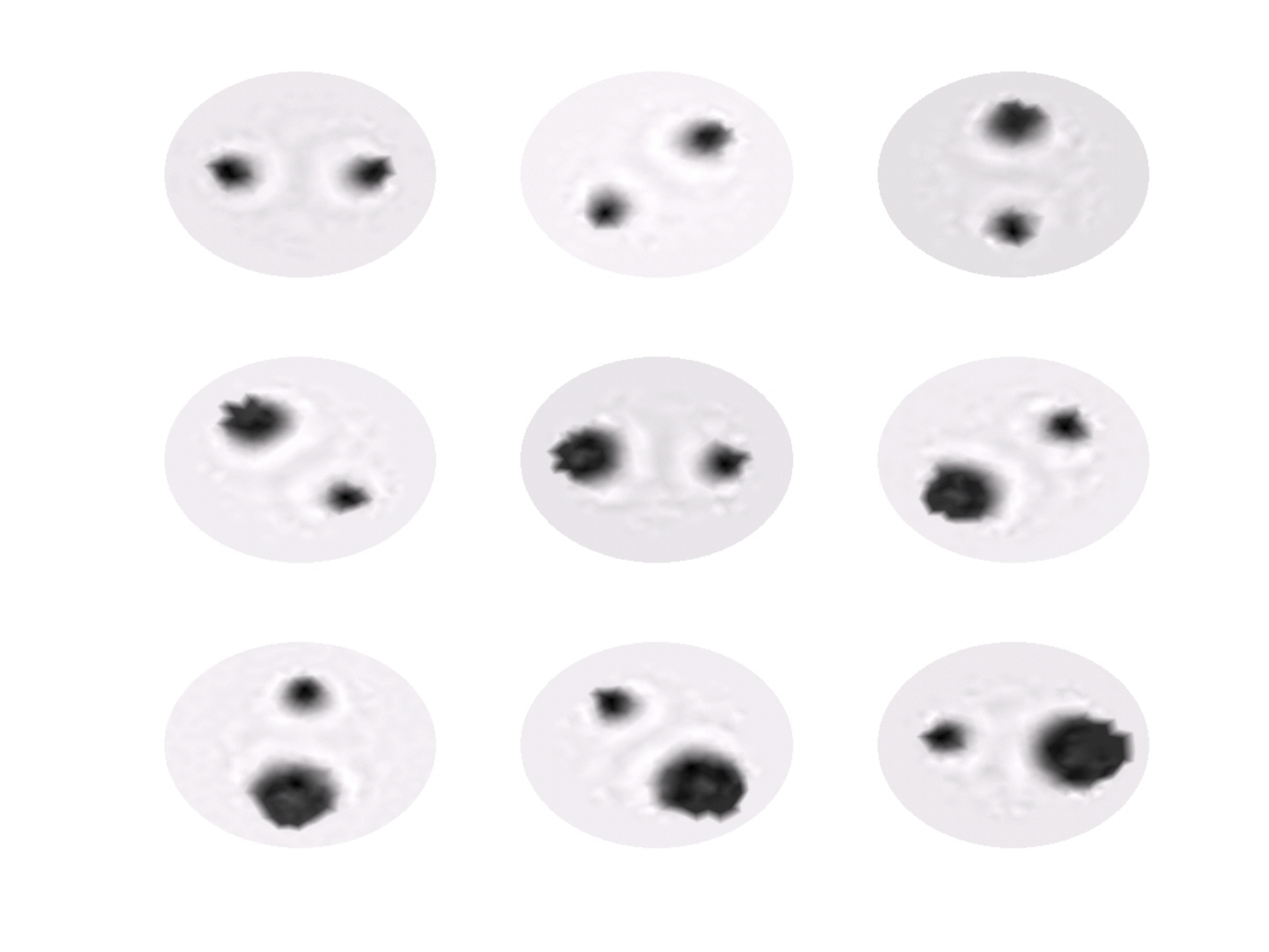}
  \label{fig:exact}
  \caption{Reconstruction results for linearized data without noise.}
\end{center}
\end{figure}

The inhomogeneity $\gamma$ was chosen as the characteristic function of
two circles one with fixed radius $r = 0.1$ the other one with increasing radius
 moving around the center of the unit ball
on an orbit of radius $0.5$. The space domain was discretized by a
triangular mesh where we used both piecewise linear finite elements for the
solution $u$ of Equation~\eqref{eq:dimp} as well as for $\gamma$. The
time-domain was discretized into $50$ uniform time-steps.
Figure~\ref{fig:exact} shows the reconstructed solution over time. The
results were computed by the discrete iterations
 defined in \eqref{itone}--\eqref{itthree}.
The location of the circles can be clearly seen from the pictures. Note that
we used an $L^2$-regularization matrix, hence the images are blurred, which
is to be expected with such a linear regularization. For a characteristic
function $\gamma$, as in our example, a bounded variation type regularization
would be suited better, but it is not clear how to incorporate such a
nonlinear regularization term into this dynamic programming framework.

\begin{figure}[ht]
\begin{center}
  \includegraphics[width=\textwidth]{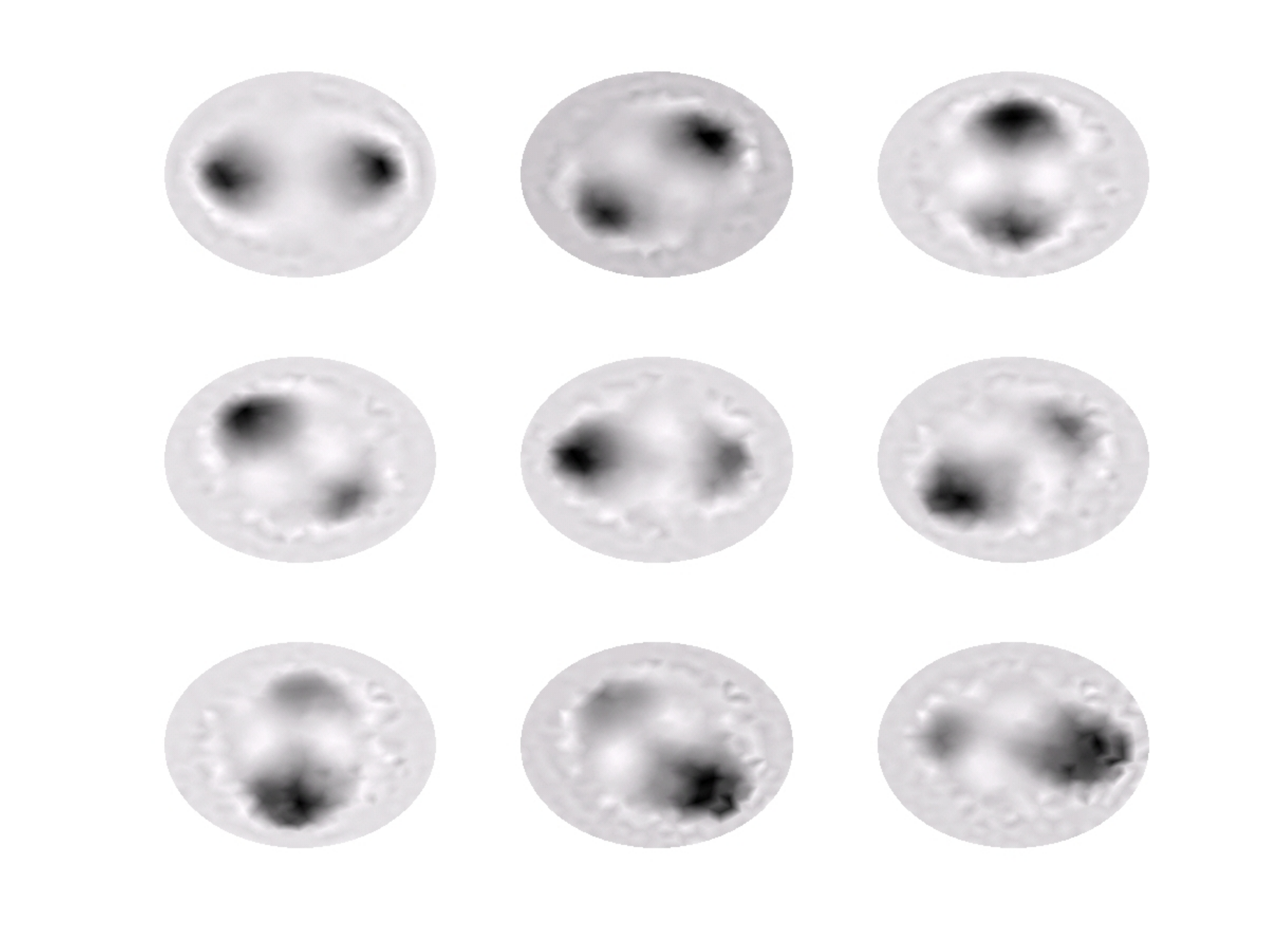}
  \label{fig:nonlinear}
  \caption{Reconstruction results for full nonlinear data with 1 \%
  random noise.}
\end{center}
\end{figure}

For the second example we used the full nonlinear data $y =
\Lambda_{1+\gamma(.,t)} - \Lambda_1 \sim F'(1) \gamma,$ and added $1 \%$
white noise to the data. But for the computation of the evolution we still
used the linearized operator $F'(1)$.
The results are shown in Figure~\ref{fig:nonlinear}. Also in this case
- even with a systematical error due to the linearization - we still get
fairly good results.

%---------------------------------------------------------------------------
% Section 5
%---------------------------------------------------------------------------
\section{Conclusions}

Each method derived in this paper require, in a first step, the solution
of an evolutionary equation (of Hamilton-Jacobi type). In a second step,
the components of the solution vector $\{ u_k \}$ are computed one at a
time. This strategy reduces significantly both the size of the systems
involved in the solution method, as well as storage requirements needed
for the numerical implementation. These points turn out to become critical
for long time measurement processes.

Some detailed considerations about complexity:
Assume that all $F(t_k)$ are discretized as $(n\times m)$ matrices.
The main effort is the matrix multiplication for the update step for
$Q_k$: In each step this requires $\mathcal{O}(n^3 + n^2 m)$ calculations.
Hence the overall complexity is of the order $\mathcal{O}(n_T(n^3 + n^2m))$
operations. If the discrete version is used, then in each step a
matrix-inversion has to be performed, which is also of the order
$\mathcal{O}(n^3)$. which leads to the same complexity as above.
In contrast, the method in \cite{SL02} requires
$\mathcal{O}((n+n_T)^3 + (n_T+m) n n_T) $. Although this is only
of cubic order in comparison to a quartic order complexity for the
dynamic programming approach, it is cubic in $n_T$. Hence if $n_T$
is large, the method proposed in this paper (which is linear in
$n_T$) will be more effective than the method in \cite{SL02}.

% The numerical results show the feasibility and the stability of our
% method. Note that the results are more smeared out at the center of the
% square, which is clear since the identification problem is less stable
% if the boundary is further away.
% 

%---------------------------------------------------------------------------
%---------------------------------------------------------------------------
\section*{Acknowledgments}
The work of S.K. is supported by Austrian Science Foundation under
grant SFB F013/F1317. % and by NSF grant Nr. DMI-0327077.
A.L. acknowledge support of CNPq under grants 305823/2003-5 and 478099/2004-5.

%---------------------------------------------------------------------------
%---------------------------------------------------------------------------

\end{document}